\documentclass[12pt]{article}
\usepackage{amsmath, amsthm, amsfonts}
\usepackage{amssymb}
\usepackage{amscd}
\usepackage{verbatim}

\begin{document}
\newcommand{\M}{{\mathcal M}}
\newcommand{\loc}{{\mathrm{loc}}}
\newcommand{\dx}{\,\mathrm{d}x}
\newcommand{\dy}{\,\mathrm{d}y}
\newcommand{\core}{C_0^{\infty}(\Omega)}
\newcommand{\sob}{W^{1,p}(\Omega)}
\newcommand{\sobloc}{W^{1,p}_{\mathrm{loc}}(\Omega)}
\newcommand{\merhav}{{\mathcal D}^{1,p}}
\newcommand{\be}{\begin{equation}}
\newcommand{\ee}{\end{equation}}
\newcommand{\mysection}[1]{\section{#1}\setcounter{equation}{0}}
\newcommand{\bea}{\begin{eqnarray}}
\newcommand{\eea}{\end{eqnarray}}
\newcommand{\bean}{\begin{eqnarray*}}
\newcommand{\eean}{\end{eqnarray*}}
\newcommand{\thkl}{\rule[-.5mm]{.3mm}{3mm}}
\newcommand{\cw}{\stackrel{\rightharpoonup}{\rightharpoonup}}
\newcommand{\id}{\operatorname{id}}
\newcommand{\supp}{\operatorname{supp}}
\newcommand{\wlim}{\mbox{ w-lim }}
\newcommand{\mymu}{{x_N^{-p_*}}}
\newcommand{\R}{{\mathbb R}}
\newcommand{\N}{{\mathbb N}}
\newcommand{\Z}{{\mathbb Z}}
\newcommand{\Q}{{\mathbb Q}}
\newcommand{\abs}[1]{\lvert#1\rvert}
\newtheorem{theorem}{Theorem}[section]
\newtheorem{corollary}[theorem]{Corollary}
\newtheorem{lemma}[theorem]{Lemma}
\newtheorem{definition}[theorem]{Definition}
\newtheorem{remark}[theorem]{Remark}
\newtheorem{proposition}[theorem]{Proposition}
\newtheorem{assertion}[theorem]{Assertion}
\newtheorem{problem}[theorem]{Problem}
\newtheorem{conjecture}[theorem]{Conjecture}
\newtheorem{question}[theorem]{Question}
\newtheorem{example}[theorem]{Example}
\newtheorem{Thm}[theorem]{Theorem}
\newtheorem{Lem}[theorem]{Lemma}
\newtheorem{Pro}[theorem]{Proposition}
\newtheorem{Def}[theorem]{Definition}
\newtheorem{Exa}[theorem]{Example}
\newtheorem{Exs}[theorem]{Examples}
\newtheorem{Rems}[theorem]{Remarks}
\newtheorem{Rem}[theorem]{Remark}

\newtheorem{Cor}[theorem]{Corollary}
\newtheorem{Conj}[theorem]{Conjecture}
\newtheorem{Prob}[theorem]{Problem}
\newtheorem{Ques}[theorem]{Question}
\newcommand{\pf}{\noindent \mbox{{\bf Proof}: }}


\renewcommand{\theequation}{\thesection.\arabic{equation}}
\catcode`@=11 \@addtoreset{equation}{section} \catcode`@=12


\title{On positive solutions of $p$-Laplacian-type equations}
\author{Yehuda Pinchover\\
 {\small Department of Mathematics}\\ {\small  Technion - Israel Institute of Technology}\\
 {\small Haifa 32000, Israel}\\
{\small pincho@techunix.technion.ac.il}\\\and Kyril Tintarev
\\{\small Department of Mathematics}\\{\small Uppsala University}\\
{\small SE-751 06 Uppsala, Sweden}\\{\small
kyril.tintarev@math.uu.se}}
 \maketitle
\newcommand{\dnorm}[1]{\thkl #1 \thkl\,}
\vspace{-1mm}
\begin{center}{\it Dedicated to Vladimir Maz'ya on the occasion of his 70th
birthday}
\end{center}

\begin{abstract}
Let $\Omega$ be a domain in $\mathbb{R}^d$, $d\geq 2$, and
$1<p<\infty$. Fix $V\in L_{\mathrm{loc}}^\infty(\Omega)$. Consider
the functional $Q$ and its G\^{a}teaux  derivative $Q^\prime$
given by
$$Q(u):=\!\!\frac{1}{p}\int_\Omega\!\!\! (|\nabla u|^p+V|u|^p)\!\dx,\;\; Q^\prime
(u):=\!-\nabla\cdot(|\nabla u|^{p-2}\nabla u)+V|u|^{p-2}\!u.$$  In
this paper we discuss several aspects of relations between
functional-analytic properties of the functional $Q$ and
properties of  positive solutions of the equation $Q^\prime
(u)=0$.
\\[2mm]
\noindent  2000  \! {\em Mathematics  Subject  Classification.}
Primary  \! 35J60; Secondary  35J20, 35J70, 49R50.\\[1mm]
\end{abstract}
\mysection{Introduction and Preliminaries}\label{sec1}
Properties of positive solutions of quasilinear elliptic
equations, and in particular of equations with the $p$-Laplacian
term in the principal part, have been extensively studied over the
recent decades, (see for example \cite{AH1,AH2,ky3,PuccS} and the
references therein). Fix $p\in(1,\infty)$, a domain
$\Omega\subseteq\R^d$ and a real valued potential $V\in
L^\infty_\loc(\Omega)$. The $p$-Laplacian equation in $\Omega$
with potential $V$ is the equation of the form
\be\label{eq}-\Delta_p(u)+V|u|^{p-2}u=0\quad \mbox{in } \Omega,\ee
where $\Delta_p(u):=\nabla\cdot(|\nabla u|^{p-2}\nabla u)$ is the
celebrated $p$-Laplacian. This equation, in the semistrong sense,
is a critical point equation for the functional \be \label{Q}
Q(u)=Q_V(u):=\frac{1}{p}\int_\Omega \left(|\nabla
u|^p+V|u|^p\right)\dx\qquad u\in\core.\ee
So, we consider solutions of \eqref{eq} in the following weak
sense.
\begin{definition}{\em A function
$v\in W^{1,p}_{\mathrm{loc}}(\Omega)$ is a {\em (weak) solution}
of the equation  \be \label{groundstate}
 Q^\prime
(u):=-\Delta_p(u)+V|u|^{p-2}u=0\quad \mbox{in }  \Omega,\ee if for
every $\varphi\in\core$
 \be \label{solution} \int_\Omega (|\nabla v|^{p-2}\nabla
v\cdot\nabla\varphi+V|v|^{p-2}v\varphi)\dx=0. \ee
We say that a real function $v\in C^1_{\mathrm{loc}}(\Omega)$ is a
{\em supersolution} (resp. {\em subsolution})  of the equation
(\ref{groundstate}) if for every nonnegative $\varphi\in\core$
 \be\label{supersolution}
\int_\Omega (|\nabla v|^{p-2}\nabla
v\cdot\nabla\varphi+V|v|^{p-2}v\varphi)\dx\geq 0 \mbox{ (resp.
}\leq 0\mbox{).} \ee
 }\end{definition}
Next, we  mention {\em local} properties of solutions of
(\ref{groundstate}) that hold in any smooth subdomain
$\Omega'\Subset\Omega$ (i.e., $\overline{\Omega'}$ is compact in
$\Omega$).

\vskip 3mm

\noindent {\bf 1. Smoothness and Harnack inequality.}  Weak
solutions of (\ref{groundstate}) admit H\"older continuous first
derivatives, and nonnegative solutions of (\ref{groundstate})
satisfy the Harnack inequality (see for example
\cite{LU,Serrin1,Serrin2,PuccS,T}).

\vskip 3mm

\noindent {\bf  2. Principal eigenvalue and eigenfunction.} For
any smooth subdomain $\Omega'\Subset\Omega$ consider the
variational problem \be \label{mu} \lambda_{1,p}(\Omega'):=\inf_{
u\in W_0^{1,p}(\Omega')}\dfrac{\int_{\Omega'}(|\nabla
u|^p+V|u|^p)\dx}{\int_{\Omega'} {|u|^p}\dx}\,.\ee
 It is well-known that for such a subdomain, (\ref{mu}) admits (up to a multiplicative
constant) a unique minimizer $\varphi$ \cite{DKN,GS}. Moreover,
$\varphi$ is a positive solution of the quasilinear eigenvalue
problem
\begin{equation}
  \begin{cases}
Q'(\varphi)=\lambda_{1,p}(\Omega')|\varphi|^{p-2}\varphi    & \text{ in } \Omega', \\
    \varphi=0 & \text{ on } \partial \Omega'.
  \end{cases}
\end{equation}
$\lambda_{1,p}(\Omega')$ and $\varphi$ are called, respectively,
the {\em principal eigenvalue and eigenfunction} of the operator
$Q'$ in $\Omega'$.

\vskip 3mm

\noindent {\bf 3. Weak and strong maximum principles.}
\begin{theorem}[{\cite{GS} (see also \cite{AH1,AH2})}]\label{thmGS}
Assume that $\Omega\subset\R^d$ is a bounded
$C^{1+\alpha}$-domain, where $0<\alpha\leq 1$. Consider a
functional $Q$ of the form \eqref{Q} with $V\in L^\infty(\Omega)$.
The following assertions are equivalent:

\begin{itemize}
 \item[(i)] $Q'$  satisfies the
maximum principle: If $u$ is a solution of the equation
$Q'(u)=f\geq 0$ in $\Omega$ with some $f\in L^\infty(\Omega)$, and
satisfies $u\geq 0$ on $\partial\Omega$, then $u$ is nonnegative
in $\Omega$.

\item[(ii)] $Q'$  satisfies the strong maximum principle: If $u$
is a solution of the equation $Q'(u)=f\gneqq 0$ in $\Omega$ with
some $f\in L^\infty(\Omega)$, and satisfies $u\geq 0$ on
$\partial\Omega$, then $u>0$ in $\Omega$.

\item[(iii)] $\lambda_{1,p}(\Omega)>0$.

\item[(iv)] For some $0\lneqq f\in L^\infty(\Omega)$ there exists a positive strict
supersolution $v$ satisfying  $Q'(v)=f$ in $\Omega$, and $v=0$ on
$\partial\Omega$.

\item[(iv')] There exists a positive strict supersolution
$v\in W^{1,p}(\Omega)\cap L^\infty(\Omega)$ satisfying
$Q'(v)=f\gneqq 0$ in $\Omega$, such that $v|_{\partial \Omega}\in
C^{1+\alpha}(\partial \Omega)$ and $f\in L^\infty(\Omega)$.

\item[(v)] For each nonnegative $f\in C^\alpha(\Omega)\cap
L^\infty(\Omega)$ there exists a unique weak nonnegative solution
of the problem $Q'(u)=f$ in $\Omega$, and $u=0$ on
$\partial\Omega$.
 \end{itemize}
 \end{theorem}


\vskip 3mm

\noindent {\bf 4. Weak comparison principle.} We recall also the
following {\em weak comparison principle} (or WCP for brevity).
\begin{theorem}[\cite{GS}]\label{CP}
Let $\Omega\subset\mathbb{R}^d$ be a bounded domain of class
$C^{1,\alpha}$, where $0< \alpha\leq 1$, and suppose that $V\in
L^\infty(\Omega)$. Assume that $\lambda_{1,p}(\Omega)>0$ and let
$u_i\in W^{1,p}(\Omega)\cap L^\infty(\Omega)$ satisfying
$Q'(u_i)\in L^\infty(\Omega)$, $u_i|_{\partial \Omega}\in
C^{1+\alpha}(\partial \Omega)$, where $i=1, 2$. Suppose further
that the following inequalities are satisfied
\begin{equation}
\left\{
\begin{array}{rclc}
 Q'(u_1)\!\!\!&\leq &\!\!\! Q'(u_2) \qquad &\mbox{ in }
\Omega,\\
Q'(u_2)\!\!\!&\geq &\!\!\! 0 \qquad &\mbox{ in }
\Omega,\\
u_{1}\!\!\!&\leq&\!\!\!u_{2} \qquad &\mbox{ on }\partial\Omega,\\
u_{2}\!\!\!&\geq&\!\!\!0 \quad&\mbox{ on } \partial \Omega.\\
\end{array}\right.
\end{equation}
Then
$$ u_{1}\leq u_{2} \qquad\mbox{ in }\Omega.$$
\end{theorem}

\vskip 3mm

\noindent {\bf 5. Strong comparison principle.}
\begin{definition}\label{defSCP}{\em we say that the {\em strong
comparison principle} (or SCP for brevity) holds true for the
functional $Q$ if the conditions of Theorem~\ref{CP} implies that
$u_{1}< u_{2}$ in $\Omega$ unless $u_1=u_2$ in $\Omega$.}
\end{definition}

\begin{remark}\label{remSCP}{\em
It is well known that the SCP holds true for $p=2$ and for
$p$-harmonic functions. For sufficient conditions for the validity
of the SCP see \cite{AS,CT,DaS,LP,PuccS,T1} and the references
therein. In \cite{CT} M.~ Cuesta and P.~Tak\'{a}\v{c} present a
counterexample where the WCP holds true but the SCP does not.}
\end{remark}

Throughout this paper we assume that \be Q(u)\ge 0 \qquad \forall
u\in \core.\ee
The following Allegretto-Piepenbrink-type theorem, links the
existence of positive solutions with the positivity of $Q$.
\begin{theorem}[{\cite[Theorem~2.3]{ky3}}]\label{pos} Consider a functional $Q$
of the form (\ref{Q}).
The following assertions are equivalent:
\begin{itemize}
 \item[(i)] The functional $Q$ is nonnegative on
$C_0^\infty(\Omega)$.
 \item[(ii)] Equation (\ref{groundstate}) admits a global positive solution.
 \item[(iii)] Equation
(\ref{groundstate}) admits a global positive supersolution.
\end{itemize}
\end{theorem}

In this paper we survey further connections between
functional-analytic properties of the functional $Q$ and
properties of its positive solutions. In particular, we review the
following topics:
\begin{itemize}
\item A representation of the nonnegative functional $Q$ as an integral
of a nonnegative Lagrangian density and the existence of a useful
equivalent nonnegative Lagrangian density with a simplified form
(Section \ref{sec:picone}).
\item The equivalence of several weak coercivity properties of $Q$.
The characterization of the non-coercive case in terms of a
Poincar\'e-type inequality, in terms of the existence of a
generalized ground state, and in terms of the variational capacity
of balls (Section \ref{sec:coercivity}).
\item  The identification of ground state as a global minimal
solution (Section \ref{sec:minimal growth}).
 \item A theorem of Liouville type connecting the
behavior of a ground state of one functional with the existence of
a ground state of another functional with a given `decaying'
subsolution (Section \ref{sec:Liouville theorems}).
\item A variational
principle that characterizes solutions of minimal growth at
infinity (Section \ref{sec:mingr1}).
\item The existence of solutions to the inhomogeneous equation
$Q'(u)=f$ in the absence of the ground state (Section
\ref{sec:inhomogeneous equation}).
\item The dependence of weak coercivity
on the potential and the domain (Section \ref{sec:crit}). In
particular, in Theorem~\ref{thm:VW}, we extend the result for
$p=2$ proved in \cite[Theorem~2.9]{ky6}.

\item Properties verified only in the linear case ($p=2$), in particular, the definition of
a natural functional space associated with the functional $Q$
(Section \ref{sec:linear case}).
\end{itemize}
\mysection{Positive Lagrangian representations}\label{sec:picone}
Let $v\in C_{\mathrm{loc}}^1(\Omega)$ be a positive solution
(resp. subsolution) of (\ref{groundstate}). Using the {\em
positive Lagrangian representation} from \cite{AH1,AH2,DS},
 we infer that for every $u\in\core$, $u\ge 0$,
 \be\label{QL} Q(u)=\int_\Omega L(u,v)\dx,\qquad \mbox{resp. }\; Q(u) \leq \int_\Omega
 L(u,v)\dx,
 \ee
where \be\label{piconeLag} L(u,v):=\frac{1}{p}\left[ |\nabla
u|^p+(p-1)\frac{u^p}{v^p}|\nabla
v|^p-p\frac{u^{p-1}}{v^{p-1}}\nabla u\cdot|\nabla v|^{p-2}\nabla
v\right].\ee It can be easily verified that $L(u,v)\ge 0$ in
$\Omega$  \cite{AH1,AH2}.

Let now $w:=u/v$, where $v$ is a positive solution  of
(\ref{groundstate}) and $u\in\core$, $u\ge 0$. Then \eqref{QL}
implies that
\be \label{QL1} Q(vw)=\frac{1}{p}\int_\Omega\left[|w\nabla
v+v\nabla w |^p-w^p|\nabla v|^p -pw^{p-1}v|\nabla v|^{p-2}\nabla
v\cdot\nabla w\right]\dx. \ee Similarly,  if $v$ is a nonnegative
subsolution of (\ref{groundstate}), then \be \label{QL1sub}
Q(vw)\leq \frac{1}{p}\int_\Omega\left[|w\nabla v+v\nabla w
|^p-w^p|\nabla v|^p -pw^{p-1}v|\nabla v|^{p-2}\nabla v\cdot\nabla
w\right]\dx. \ee
Therefore, a nonnegative functional $Q$ can be represented as the
integral of a nonnegative Lagrangian $L$. In spite of the
nonnegativity of the in \eqref{QL} and \eqref{QL1}, the expression
\eqref{piconeLag} of $L$ contains an indefinite term which poses
obvious difficulties for extending the domain of the functional to
more general weakly differentiable functions. The next proposition
shows that $Q$ admits a two-sided estimate by a simplified
Lagrangian containing only nonnegative terms. We call the
functional associated with this simplified Lagrangian the {\em
simplified energy}.

Let $f$ and $g$ be two nonnegative functions. We denote $f\asymp
g$ if there exists a positive constant $C$ such that $C^{-1}g\leq
f \leq Cg$.
\begin{proposition}[{\cite[Lemma~2.2]{aky}}]
\label{prop:superPicone}  Let $v\in C_{\mathrm{loc}}^1(\Omega)$ be
a positive solution of (\ref{groundstate}). Then
\bea \label{p<2} Q(vw)\asymp  \int_\Omega\!\! v^2 |\nabla
w|^2\left(w|\nabla v|+v|\nabla w|\right)^{p-2}\dx \quad \forall
w\in C^1_0(\Omega), w\geq 0. \eea
In particular, for $p\ge 2$, we have
\bea \label{p>2} Q(vw)\asymp  \int_\Omega \!\!\left(v^p|\nabla
w|^p+v^2|\nabla v|^{p-2} w^{p-2}|\nabla w|^2\right)\!\!\dx \quad
\forall w\in C^1_0(\Omega), w\geq 0. \eea

If $v$ is only a nonnegative subsolution of (\ref{groundstate}),
then for $1<p<\infty$ we have
\bean \label{p<2sub} Q(vw)\!\leq\! C\!\!  \int_{\Omega\cap\{v>
0\}}\!\!\!\! v^2 |\nabla w|^2\left(w|\nabla v|+v|\nabla
w|\right)^{p-2}\!\! \dx \quad \forall w\in C^1_0(\Omega),  w\geq
0. \eean
 In particular, for $p\ge 2$ we
have
\bean \label{p>2sub}
 Q(vw)\leq C  \!\int_\Omega \!\!\!\left(v^p|\nabla
w|^p+v^2|\nabla v|^{p-2} w^{p-2}|\nabla w|^2\right)\!\! \dx \quad
\forall w\in C^1_0(\Omega), w\geq 0. \eean
\end{proposition}
%
\begin{remark}{\em  It is shown in \cite{aky} that for $p>2$
neither of the terms in the simplified energy \eqref{p>2} is
dominated by the other, so that \eqref{p>2} cannot be further
simplified.}
\end{remark}
\mysection{Coercivity and ground state} \label{sec:coercivity} It
is well known (see \cite{M86}) that for a nonnegative Schr\"odinger operator $P$ we
have the following dichotomy: either there exists a strictly
positive potential $W$ such that the Schr\"odinger operator $P-W$
is nonnegative, or $P$ admits a unique (generalized) ground state.
It turns out that this statement is also true for nonnegative
functionals of the form \eqref{Q}.

\begin{definition}{\em  Let $Q$ be a nonnegative functional on
$\core$ of the form \eqref{Q}. We say that a sequence
$\{u_k\}\subset\core$ of nonnegative functions is a {\em null
sequence} of the functional $Q$ in $\Omega$, if there exists an
open set $B\Subset\Omega$ such that $\int_B|u_k|^p\dx=1$, and \be
\lim_{k\to\infty}Q(u_k)=\lim_{k\to\infty}\int_\Omega (|\nabla
u_k|^p+V|u_k|^p)\dx=0.\ee
 We say that a positive function $v\in
C^1_{\mathrm{loc}}(\Omega)$ is a {\em ground state} of the
functional $Q$ in $\Omega$ if $v$ is an
$L^p_{\mathrm{loc}}(\Omega)$ limit of a null sequence of $Q$. If
$Q\ge 0$, and $Q$ admits a ground state in $\Omega$, we say that
$Q$ is {\em critical} in $\Omega$. }
\end{definition}
\begin{remark}\label{remc1}{\em
The requirement that $\{u_k\}\subset \core$, can clearly be
weakened by assuming only that $\{u_k\}\subset W^{1,p}_0(\Omega)$.
Also, as it follows from Theorem~\ref{thmky3}, the requirement
$\int_B|u_k|^p\dx=1$ can be replaced by $\int_B|u_k|^p\dx\asymp 1$
or by $\int_B u_k\dx\asymp 1$.}
\end{remark}
The following statements are based on rephrased statements of
\cite[Theorem~1.6]{ky3}, \cite[Theorem~4.3]{ky5} and
\cite[Proposition~3.1]{TakTin} (cf. \cite{M86,ky2} for the case
$p=2$).
\begin{theorem}
\label{thmgs} Suppose that the functional $Q$ is nonnegative on
$\core$.
\begin{enumerate}
  \item Any ground state $v$ is a positive solution of
(\ref{groundstate}).
  \item $Q$ admits a ground state $v$ if and only if
(\ref{groundstate}) admits a unique positive supersolution.
    \item  $Q$ is critical in
$\Omega$ if and only if $Q$ admits a null sequence that converges
locally uniformly in $\Omega$.
  \item If $Q$ admits a ground state $v$, then the following Poincar\'e type
inequality holds: There exists a positive continuous function $W$
in $\Omega$, such that for every $\psi\in C_0^\infty(\Omega)$
satisfying $\int \psi v \,\mathrm{d}x \neq 0$ there exists a
constant $C>0$ such that the following inequality holds:
 \be\label{Poinc}
 Q(u)+C\left|\int_{\Omega} \psi
 u\,\mathrm{d}x\right|^p \geq
 C^{-1}\int_{\Omega} W\left(|\nabla u|^p+|u|^p\right)\,\mathrm{d}x
  \qquad
 \forall u\in C_0^\infty(\Omega).\ee
\end{enumerate}
\end{theorem}
The following theorem slightly extends \cite[Theorem~1.6]{ky3} in
the spirit of \cite[Proposition~3.1]{TakTin}.
\begin{theorem}
\label{thmky3}  Suppose that the functional $Q$ is nonnegative on
$\core$. The following statements are equivalent.
\begin{itemize}
\item[(a)] $Q$ does not admit a ground state in $\Omega$.

\item[(b)] There exists a continuous function $W>0$  in $\Omega$
such that
 \be \label{gap_p}
Q(u)\ge \int_\Omega W(x)|u(x)|^p\dx\qquad \forall u\in\core.\ee
\item[(c)]
There exists a continuous function $W>0$  in $\Omega$ such that
\be \label{gap_max} Q(u)\ge \int_\Omega W(x)\left(|\nabla
u(x)|^p+|u(x)|^p\right)\dx \qquad \forall u\in\core. \ee
\item[(d)] There exists an open set $B\Subset\Omega$ and $C_B>0$
such that \be Q(u)\ge C_B\left|\int_B
u(x)\dx\right|^p\label{gap_min}\qquad \forall u\in\core.\ee
\end{itemize}
Suppose further that $d>p$. Then $Q$ does not admit a ground state
in $\Omega$ if and only if there exists a continuous function
$W>0$  in $\Omega$ such that \be \label{gap_p*}
Q(u)\ge\left(\int_\Omega W(x)|u(x)|^{p^*}\dx\right)^{p/p^*} \qquad
\forall u\in\core, \ee where $p^*={pd}/{(d-p)}$ is the critical
Sobolev exponent.
\end{theorem}

\begin{definition} {\em
A nonnegative functional $Q$ on $\core$ of the form \eqref{Q}
which is not critical is said to be {\em subcritical} (or {\em
weakly coercive}) in $\Omega$.}
\end{definition}
\begin{example}\label{ex1}{\em Consider the functional
$Q(u):=\int_{\mathbb{R}^d}|\nabla u|^p\dx$. It follows from
\cite[Theorem~2]{MP} that if $d\leq p$, then $Q$ admits a ground
state $\varphi=\mathrm{constant}$ in $\mathbb{R}^d$. On the other
hand, if $d>p$, then
$$u(x):=\left[1+|x|^{p/(p-1)}\right]^{(p-d)/p}, \qquad v(x):=\mathrm{constant}$$
are two positive supersolutions of the equation $-\Delta_pu=0$ in
$\mathbb{R}^d$. Therefore, $Q$ is weakly coercive
in $\mathbb{R}^d$.
 }\end{example}

\begin{example}
\label{ex2} {\em Let $d>1$, $d\neq p$, and
$\Omega:=\mathbb{R}^d\setminus \{0\}$ be the punctured space. The
following celebrated Hardy's inequality holds true:
\begin{equation}\label{Hardy_fun}
  Q_\lambda(u):=\int_\Omega \left(|\nabla u|^p-
  \lambda\frac{|u|^p}{|x|^p}\right)\dx\ge 0
 \qquad u\in \core,
\end{equation}
whenever $\lambda\le c^*_{p,d}:=|{(p-d)}/{p}|^p$. Clearly, if
$\lambda< c^*_{p,d}\,$, then $Q_\lambda(u)$ is weakly coercive. On
the other hand, the proof of Theorem 1.3 in \cite{PS} shows that
$Q_\lambda$ with $\lambda=c^*_{p,d}$ admits a null sequence. It
can be easily checked that the function $v(r):=|r|^{(p-d)/p}$ is a
positive solution of the corresponding radial equation:
$$-|v'|^{p-2}\left[(p-1)v''+\frac{d-1}{r}v'\right]-c^*_{p,d}\frac{|v|^{p-2}v}{r^p}
=0\qquad r\in (0,\infty).$$ Therefore, $\varphi(x):=|x|^{(p-d)/p}$
is the ground state of the equation
\begin{equation}\label{eqHardy}
-\Delta_p u- c^*_{p,d}\frac{|u|^{p-2}u}{|x|^{p}}=0 \qquad \mbox{in
} \Omega.
\end{equation}
Note that $\varphi\not \in W^{1,p}_{\mathrm{loc}}(\R^d)$ for
$p\neq d$. In particular, $\varphi$ is not a positive
supersolution of the equation $\Delta_pu=0$ in $\R^d$. }
\end{example}

In \cite{Tr1,Tr2} Troyanov has established a relationship between
the $p$-capacity of closed balls in  a Riemannian manifold
$\mathcal{M}$ and the $p$-parabolicity of $\mathcal{M}$ with
respect the $p$-Laplacian. We extend his definition and result to
our case.
\begin{definition} {\em
Suppose that the functional $Q$ is nonnegative on $\core$. Let
$K\Subset \Omega$ be a compact set. The {\em $Q$-capacity} of $K$
in $\Omega$ is defined by
$$\mathrm{Cap}_Q(K,\Omega):=\inf\{Q(u)\mid u\in \core,\; u\geq 1 \;\mbox{ on } K\}.$$
}
\end{definition}
\begin{corollary} Suppose that the functional $Q$ is nonnegative on
$\core$. Then $Q$ is critical in $\Omega$ if and only if the
$Q$-capacity of each closed ball in $\Omega$ is zero.
 \end{corollary}

\mysection{Ground states and minimal growth at
infinity}\label{sec:minimal growth}
\begin{definition} {\em
Let $K_0$ be a compact set in $\Omega$.  A positive solution $u$
of the equation $Q'(u)=0$ in $\Omega\setminus K_0$ is said to be a
{\em positive solution of minimal growth in a neighborhood of
infinity in} $\Omega$ (or $u\in\M_{\Omega,K_0}$ for brevity) if
for any compact set $K$ in $\Omega$, with a smooth boundary, such
that $K_0 \Subset \mathrm{int}(K)$, and any positive supersolution
$v\in C((\Omega\setminus K)\cup
\partial K)$ of the equation $Q'(u)=0$ in $\Omega\setminus K$,
the inequality $u\le v$ on $\partial K$ implies that $u\le v$ in
$\Omega\setminus K$.

A (global) positive solution $u$ of the equation $Q'(u)=0$ in
$\Omega$, which has minimal growth in a neighborhood of infinity
in $\Omega$ (i.e. $u\in\M_{\Omega,\emptyset}$) is called a {\em
global minimal solution of the equation $Q'(u)=0$ in $\Omega$}.}
\end{definition}
\begin{theorem}[{\cite[Theorem~5.1]{ky5}, cf. \cite{ky3}}]\label{thmmingr1}
Suppose that $1<p<\infty$, and $Q$ is nonnegative on $\core$. Then
for any $x_0\in \Omega$ the equation $Q'(u)=0$ has a positive
solution $u\in\M_{\Omega,\{x_0\}}$.
\end{theorem}
We have the following connection between the existence of a global
minimal solutions and weak coercivity.

\begin{theorem}[{\cite[Theorem~5.2]{ky5}, cf. \cite{ky3}}]\label{thmmingr13}
Let $1<p<\infty$, and assume that $Q$ is nonnegative on $\core$.
Then $Q$ is subcritical in $\Omega$ if and only if the equation
$Q'(u)=0$ in $\Omega$ does not admit a global minimal solution of
the equation $Q'(u)=0$ in $\Omega$. In particular, $u$ is ground
state of the equation $Q'(u)=0$ in $\Omega$ if and only $u$ is a
global minimal solution of this equation.
\end{theorem}
Consider a positive solution $u$ of the equation $Q'(u)=0$ in a
punctured neighborhood of $x_0$ which has a nonremovable
singularity at $x_0\in \mathbb{R}^d$. Without loss of generality
we may assume that $x_0=0$. If  $1<p\leq d$, then the behavior of
$u$ near an isolated singularity is well understood.  Indeed, due
to a result of L.~V\'{e}ron  (see \cite[Lemma~5.1]{ky3}), we have
that
 \begin{equation}\label{nonremovasymp}
  u(x)\sim\begin{cases}
    \abs{x}^{\alpha(d,p)} & p<d, \\
     -\log \abs{x} & p=d,
  \end{cases} \qquad \mbox{ as } x\to 0,
\end{equation}
where $\alpha(d,p):=(p-d)/(p-1)$, and $f\sim g$ means that $$
\lim_{x\to 0}\frac{f(x)}{g(x)}= C$$ for some positive constant
$C$. In particular,  $\lim_{x\to 0} u(x)=\infty$.

Assume now that $p>d$. A general question is whether in this case,
any positive solution of the equation $Q'(u)=0$ in a punctured
ball centered at $x_0$ can be continuously extended at $x_0$ (see
\cite{Man} for partial results).

Under the assumption that  $u\asymp 1$ near the isolated point the
answer is given by Lemma~5.3 in \cite{ky5}:

\begin{lemma}\label{lemveron}
Assume that $p> d$, and let  $v$ be a positive solution of the
equation $Q'(u)=0$ in a punctured neighborhood of $x_0$ satisfying
$u\asymp 1$ near $x_0$.  Then $u$ can be continuously extended at
$x_0$.
\end{lemma}




The following statement combines the second part of
\cite[Theorem~5.4]{ky3}, where the case $1<p\leq d$ is considered
with \cite[Theorem~5.3]{ky5} which deals with the case $p>d$.

\begin{theorem}\label{cor_nonremove} Let $x_0\in \Omega$, and let $u\in\M_{\Omega,\{x_0\}}$.
Then $Q$ is subcritical in $\Omega$ if and only if $u$ has a
nonremovable singularity at $x_0$.
\end{theorem}

\mysection{Liouville theorems} \label{sec:Liouville theorems} In
\cite{aky} we use some of the positivity properties of the
nonnegative functional $Q$ discussed in the previous sections to
prove a Liouville comparison principle for equations $Q'(u)=0$ in
$\Omega$. (see Theorem~\ref{thm:p=2} for the case $p=2$).
\begin{theorem}\label{thm:main} Let $\Omega$ be a domain in $\R^d$, $d\geq 1$, and let
$p\in(1,\infty)$. For $j=0,1$, let $V_j\in L^\infty_\loc(\Omega)$,
and  let
$$ Q_j(u):=\int_\Omega\left(|\nabla
u(x)|^p+V_j(x)|u(x)|^p\right)\dx \qquad u\in\core. $$
\par
 Assume that the following assumptions hold true.
\begin{itemize}
\item[(i)] The functional  $Q_1$ admits a ground state
 $\varphi$ in $\Omega$.

\item[(ii)]  $Q_0\geq 0$ on $\core$, and the equation $Q_0'(u)=0$ in $\Omega$ admits a
subsolution $\psi\in W^{1,p}_{\mathrm{loc}}(\Omega)$ satisfying
$\psi_+\neq 0$, where $\psi_+(x):=\max\{0, \psi(x)\}$.

\item[(iii)] The following inequality holds in $\Omega$
\begin{equation} \label{ineq:phi-psi} \psi_+\le C\varphi, \end{equation} where $C>0$ is a
positive constant.

\item[(iv)] The following inequality holds in $\Omega$
\begin{equation} \label{ineq:gradphi-gradpsi}
    |\nabla\psi_+|^{p-2}\le C|\nabla
\varphi|^{p-2},
 \end{equation} where $C>0$ is a positive constant.
\end{itemize}
Then the functional $Q_0$ admits a ground state in $\Omega$, and
$\psi$ is the ground state. In particular, $\psi$ is (up to a
multiplicative constant) the unique positive supersolution of the
equation $Q'_0(u)=0$ in $\Omega$.
\end{theorem}
\begin{remark}{\em
Condition \eqref{ineq:gradphi-gradpsi} is redundant for $p=2$. For
$p\neq 2$ it is equivalent to the assumption that the following
inequality holds in $\Omega$:
\begin{equation} \label{ineq:gradphi-gradpsi1}
 \begin{cases}
    |\nabla\psi_+|\le C|\nabla
\varphi| & \text{ if }\; p>2, \\
    |\nabla\psi_+|\ge C|\nabla
\varphi| & \text{ if }\; p<2,
  \end{cases}
 \end{equation} where $C>0$ is a positive constant.}
\end{remark}
\begin{remark}\label{remgradcond}{\em
This theorem holds if, in addition to  \eqref{ineq:phi-psi}, one
assumes instead of $|\nabla\psi_+|^{p-2}\le C|\nabla
\varphi|^{p-2}$  in $\Omega$ (see \eqref{ineq:gradphi-gradpsi}),
that the following inequality holds true in $\Omega$
\begin{equation} \label{ineq:gradphi-gradpsi-a}
    \psi_+^2|\nabla\psi_+|^{p-2}\le C\varphi^2|\nabla
\varphi|^{p-2},
 \end{equation} where $C>0$ is a positive constant.
 }
 \end{remark}
\begin{remark}\label{remNogradcond} {\em Suppose that $1<p<2$, and assume that the ground state
$\varphi>0$ of the functional $Q_1$ is such that $w=\mathbf{1}$ is
a ground state of the functional \be \label{E-1}
E_1^{\varphi}(w)=\int_\Omega \varphi^p|\nabla w|^p \dx , \ee that
is, there is a sequence $\{w_k\}\subset C_0^\infty(\Omega)$ of
nonnegative functions satisfying $E_1^{\varphi}(w_k)\to 0$, and
$\int_B |w_k|^p=1$ for a fixed $B\Subset\Omega$ (this implies that
$w_k\to\mathbf{1}$ in $L^p_{\mathrm{loc}}(\Omega)$). In this case,
the conclusion of Theorem~\ref{thm:main} holds if there is a
nonnegative subsolution $\psi_{+}$ of $Q'_0(u)=0$ satisfying
\eqref{ineq:phi-psi} alone, without an assumption on the gradients
(like \eqref{ineq:gradphi-gradpsi} or
\eqref{ineq:gradphi-gradpsi-a}).}
\end{remark}

\begin{remark}{\em
Condition \eqref{ineq:gradphi-gradpsi} is  essential when $p>2$,
and presumably also when $p<2$. When $p>2$, $\Omega=\R^d$ and $V$
is radially symmetric, Proposition~4.2 in \cite{aky} shows
that the simplified energy functional is not equivalent to either
of its two terms that lead to conditions \eqref{ineq:phi-psi} and
\eqref{ineq:gradphi-gradpsi}, respectively.
 }
 \end{remark}

\begin{example}\label{ex12}{\em
Assume that $1\leq d\leq p\leq 2$, $p>1$,   $\Omega
=\mathbb{R}^d$, and consider the functional
$Q_1(u):=\int_{\mathbb{R}^d} |\nabla u|^p\dx$. By
Example~\ref{ex1}, the functional  $Q_1$ admits a ground state
$\varphi=\mathrm{constant}$ in $\mathbb{R}^d$.

Let $Q_0$ be a functional of the form \eqref{Q} satisfying
$Q_0\geq 0$ on $C_0^\infty(\mathbb{R}^d)$. Let   $\psi\in
W^{1,p}_{\mathrm{loc}}(\mathbb{R}^d)$, $\psi_+\neq 0$ be a
subsolution of the equation $Q_0'(u)=0$ in $\mathbb{R}^d$, such
that $\psi_+\in L^\infty(\mathbb{R}^d)$. It follows from
Theorem~\ref{thm:main} that $\psi$ is the ground state of $Q_0$ in
$\mathbb{R}^d$. In particular, $\psi$ is (up to a multiplicative
constant) the unique positive supersolution and unique bounded
solution of the equation $Q'_0(u)=0$ in $\mathbb{R}^d$. Note that
there is no assumption on the behavior of the potential $V_0$ at
infinity. This result generalizes some striking Liouville theorems
for Schr\"odinger operators on $\mathbb{R}^d$ that hold for
$d=1,2$ and $p=2$ (see \cite[theorems~1.4--1.6]{Pl}).
 }
 \end{example}
\begin{example}\label{ex12a}{\em  Let $1<p<\infty$, $d>1$, $d\neq p$,
and $\Omega:=\mathbb{R}^d\setminus \{0\}$ be the punctured space.
Let $Q_0$ be a functional of the form \eqref{Q} satisfying
$Q_0\geq 0$ on $\core$. Let   $\psi\in
W^{1,p}_{\mathrm{loc}}(\Omega)$, $\psi_+\neq 0$ be a subsolution
of the equation $Q_0'(u)=0$ in $\Omega$, satisfying
\be\label{HardyPsi} \psi_+(x)\le C|x|^{(p-d)/p}\qquad
x\in\Omega.\ee When $p>2$, we require in addition that the
following inequality is satisfied \be \label{HardyGrad}
\psi_+(x)^2|\nabla\psi_+(x)|^{p-2}\le C|x|^{2-d}\qquad x\in\Omega.
\ee It follows from Theorem~\ref{thm:main},
Remark~\ref{remgradcond}, Remark~\ref{remNogradcond} and
Example~\ref{ex2} that $\psi$ is the ground state of $Q_0$ in
$\Omega$. Let $\varphi$ be the ground state of the Hardy
functional. The reason that \eqref{HardyGrad} is stated only for
$p>2$ hinges on the property that for $1<p<2$ the functional \be
\label{HardyGS} E_1^\varphi(w)=\int_{\Omega}|x|^{p-d}|\nabla
w|^p\dx \ee admits a ground state $\mathbf{1}$, so
Remark~\ref{remNogradcond} applies.}
 \end{example}
Next, we present a family of functionals $Q_0$ for which the
conditions of Example~\ref{ex12a} are satisfied.
\begin{example}
{\em Let $d\geq 2$, $1<p<d$, $\alpha\ge 0$, and
$\Omega:=\mathbb{R}^d\setminus \{0\}$. Let
$$
W_\alpha(x):=-\left(\frac{d-p}{p}\right)^p \;\dfrac{{\alpha
dp}/{(d-p)}+|x|^\frac{p}{p-1}}{\left(\alpha+|x|^\frac{p}{p-1}\right)^p}\;.
$$
Note that if $\alpha=0$ this is the Hardy potential as in
Example~\ref{ex2}. If $Q_0$ is the functional \eqref{Q} with the
potential $V_0:=W_\alpha$, then
$$
\psi_\alpha(x):=\left(\alpha+|x|^\frac{p}{p-1}\right)^{-\frac{(d-p)(p-1)}{p^2}}
$$
is a solution of $Q_0'(u)=0$ in $\Omega$, and therefore $Q_0\ge 0$
on $\core$. Moreover, one can use Example~\ref{ex12a} to show that
$\psi_\alpha$ is a ground state of $Q_0$. Note first that
$\psi=\psi_\alpha$ satisfies \eqref{HardyPsi}. If
$\frac{d}{d-1}<p<d$, then $\psi_\alpha$ satisfies also
\eqref{HardyGrad} and therefore, it is a ground state in this
case. In the remaining case $p\le\frac{d}{d-1}\leq 2$,
Example~\ref{ex2} concludes that $\psi_\alpha$ is a ground state
from the property of the functional \eqref{HardyGS}. }
\end{example}

\mysection{Variational principle for solutions of minimal growth
and comparison principle} \label{sec:mingr1} The aim of this
section is to represent positive solutions of minimal growth in a
neighborhood of infinity in $\Omega$ as a limit of a modified null
sequence.
\begin{theorem}[{\cite[Theorem~7.1]{ky5}}]\label{thmmin_null}
Suppose that $1<p<\infty$, and let $Q_V$ be nonnegative on
$\core$. Let $\Omega_1\Subset \Omega$ be an open set, and let
$u\in C(\Omega\setminus \Omega_1)$ be a positive solution of the
equation $Q_V'(u)=0$ in $\Omega\setminus \overline{\Omega_1}$
satisfying $|\nabla u|\neq 0$ in $\Omega\setminus
\overline{\Omega_1}$.
\par
Then $u\in\M_{\Omega,\overline{\Omega_1}}$ if for every smooth
open set $\Omega_2$ satisfying $\Omega_1 \Subset \Omega_2\Subset
\Omega$, and an open set $B\Subset(\Omega\setminus
\overline{\Omega_2})$ there exists a sequence
$\{u_k\}\subset\core$,  $u_k\ge 0$, such that for all $k\in\N$,
$\int_B|u_k|^p\dx=1$, and \be
\lim_{k\to\infty}\int_{\Omega\setminus \overline{\Omega}_2}
L(u_k,u)\dx =0,\ee where $L$ is the Lagrangian given by
\eqref{piconeLag}.
\end{theorem}
\begin{conjecture}\label{remNecessity1}{\em
We conjecture that for $p\neq 2$ a positive global solution of the
equation $Q'_V(u)=0$ in $\Omega$ satisfying  $u\in
\M_{\Omega,\overline{\Omega_1}}$ for some smooth open set
$\Omega_1 \Subset \Omega$  is a global minimal solution.
 }\end{conjecture}
\begin{remark}\label{remNecessity}{\em
The validity of Conjecture~\ref{remNecessity1} seems to be related
to the SCP. We note that if Conjecture~\ref{remNecessity1} holds
true, then the condition of Theorem~\ref{thmmin_null} is also
necessary (cf. Section~\ref{ssec:mingr2}).
 }\end{remark}
Finally, we formulate a sub-supersolution comparison principle for
our singular elliptic equation.
\begin{theorem}[{Comparison Principle \cite{ky5}}]
\label{thm_comparison} Assume that the functional $Q_V$ is
nonnegative on $\core$. Fix smooth open sets $\Omega_1\Subset
\Omega_2\Subset \Omega$. Let $u,v\in
W^{1,p}_{\mathrm{loc}}(\Omega\setminus \Omega_1)\cap
C(\Omega\setminus \Omega_1)$ be, respectively, a positive
subsolution and a supersolution of the equation $Q'_V(w)=0$ in
$\Omega\setminus \overline{\Omega_1}$ such that $u\leq v$ on
$\partial \Omega_2$.
\par
Assume further that $Q'_V(u)\in
L^\infty_{\mathrm{loc}}(\Omega\setminus \Omega_1)$,  $|\nabla
u|\neq 0$ in $\Omega\setminus \overline{\Omega_1}$, and that there
exist an open set $B\Subset(\Omega\setminus \overline{\Omega_2})$
and a sequence $\{u_k\}\subset\core$, $u_k\ge 0$, such that  \be
\int_B|u_k|^p\dx=1 \quad \forall k\geq 1, \mbox{ and }\quad
\lim_{k\to\infty}\int_{\Omega\setminus \overline{\Omega}_1}
L(u_k,u)\dx =0.\ee Then $u\leq v$ on $\Omega\setminus \Omega_2$.
\end{theorem}
\mysection{Solvability of nonhomogeneous
equation}\label{sec:inhomogeneous equation}
In this section we
discuss some results of \cite{TakTin} concerning the solvability
of the nonhomogeneous equation
 \be \label{Qf} Q'_V(u)=f \qquad \text{ in
}\Omega, \ee where $Q_V$ is the nonnegative functional \eqref{Q}.
In some cases, e.g. $V\ge 0$ or $p=2$, the nonnegativity of $Q_V$
implies that $Q_V$ is convex. In general, however, $Q_V$ might be
nonconvex.  For $p>2$, see the elementary one-dimensional example
at the end of \cite{dPEM}, and also the proof of
\cite[Theorem~7]{GS}. For $p<2$, see Example~2 in \cite{FHTdT}.

If $Q_V$ is convex and weakly coercive, then \eqref{Qf} can be
easily solved by defining a Banach space $X$ as a completion of
$\core$ with respect to the norm $Q_V(\cdot)^{1/p}$ (see the
discussion of the analogous space for $p=2$ in
Section~\ref{sec:linear case} below). Such space is continuously
imbedded into $W^{1,p}_{\mathrm{loc}}(\Omega)$ by \eqref{gap_max}.
It follows that for every $f\in X^*$ the functional
$$
u\mapsto Q_V(u)-\langle f,u\rangle \qquad u\in X,
$$
has a minimum that solves \eqref{Qf}.

Note that the requirement of
weak coercivity cannot be removed. Indeed, if $p=2$, $\Omega$ is
smooth and bounded, and $V=0$, then the corresponding ground state
$\varphi$ is the first eigenfunction of the Dirichlet Laplacian
with an eigenvalue $\lambda_0$, and there is no solution $u\in
W^{1,2}_0(\Omega)$ to the equation \be \label{gslapl}
(-\Delta-\lambda_0)u=f \qquad \text{ in }\Omega \ee  unless
$\int_\Omega \varphi(x)f(x)\dx=0$.

In order to address the nonconvex case, we use the following setup
from convex analysis (see \cite[Chapt.~I]{Eke-Temam} for details.)
The {\em polar\/} (or {\em conjugate\/}) functional to ${Q}_V$ is
defined by
\begin{equation}
 {Q}_V^{*}(f):=
    \sup_{ u\in \core }
    \left[ \langle u,f\rangle - {Q}_V(u) \right] \qquad f\in \mathcal{D}'(\Omega).
\label{def:Q_V^*}
\end{equation}
Notice that ${Q}_V$ is positively homogeneous of degree $p$, and
consequently, ${Q}_V^{*}$ is positively homogeneous functional of
degree $p'$, where $p'= p/(p-1)$. The ({\em effective\/} or {\em
natural\/}) {\em domain\/} $X^{*}$ of ${Q}_V^{*}$ is defined
naturally by
\begin{equation}
  X^{*} =
  \{ f\in \mathcal{D}'(\Omega)\colon {Q}_V^{*}(f) < \infty \} .
\label{X^*=dom_Q_V^*}
\end{equation}
The definition of ${Q}_V^{*}(f)$ in \eqref{def:Q_V^*} yields
immediately the well known Fenchel\--Young inequality
\begin{equation}
  | \langle u,f\rangle |\leq
    {Q}_V(u) + {Q}_V^{*}(f),
\label{ineq:Young}
\end{equation}
and equivalently, the H\"older inequality
\begin{equation}
  | \langle u,f\rangle |\leq
    \left( p\, {Q}_V(u) \right)^{1/p}\,
    \left( p'\,{Q}_V^{*}(f) \right)^{1/p'}.
\label{ineq:Hoelder}
\end{equation}
One can easily verify that $X^{*}$ is a linear subspace of
$\mathcal{D}'(\Omega)$ and
\begin{equation}
  \| f\|_{*}:= \left(p'{Q}_V^{*}(f)\right)^{1/p'}
\label{def:norm*}
\end{equation}
defines a norm on $X^{*}$. In particular, $\| f\|_{*} = 0$ implies
$f=0$, as a consequence of \eqref{ineq:Hoelder} combined with the
separation property of the duality between $\core$ and
$\mathcal{D}'(\Omega)$.

From \eqref{gap_max}
one immediately deduces that $X^{*}$ contains
$\left(W_0^{1,p}(\Omega; W)\right)'$
and that the corresponding embedding
\begin{equation*}
 \left(W_0^{1,p}(\Omega; W)\right)' \hookrightarrow X^{*}
\end{equation*}
is continuous and dense.
The density follows from
\begin{equation}
  \core \subset
\left(W_0^{1,p}(\Omega; W)\right)'
  \subset X^{*}\subset \mathcal{D}'(\Omega) .
\label{emb:L^p'->X^*}
\end{equation}
Therefore, denoting by $X^{**}$ the (strong) dual space of $X^{*}$
with respect to the duality $\langle\,\cdot\, ,\,\cdot\,\rangle$,
we observe that $X^{**}$ is continuously embedded into
$W^{1,p}(\Omega; W)$ and that $\core$ is weak\--star dense in
$X^{**}$. It is noteworthy that the separability of $X^{*}$ in the
norm topology implies that the weak\--star topology on any bounded
subset of $X^{**}$ is metrizable (Rudin \cite[Theorem 3.16,
p.~70]{Rudin-funct}). Now consider the {\em bipolar\/} (or {\em
second conjugate\/}) functional to ${Q}_V$ defined by
\begin{equation}
  {Q}_V^{**}(u):=
    \sup_{f\in X^{*}}
    \left[ \langle u,f\rangle - {Q}_V^{*}(f) \right] \qquad
    u\in X^{**}.
\label{def:Q_V^**}
\end{equation}
From \eqref{ineq:Young} it is evident that
\begin{equation}
  0\leq {Q}_V^{**}(u) \leq {Q}_V(u)
    \quad\mbox{ for every } u\in \core .
\label{ineq:Q_V^**<Q_V}
\end{equation}
Moreover, in analogy with the norm $\|\cdot\|_{*}$ on $X^{*}$ (see
\eqref{def:norm*}), the dual norm $\|\cdot\|_{**}$ on $X^{**}$ is
given by
\begin{equation}
  \| u\|_{**} = \left( p\, {Q}_V^{**}(u) \right)^{1/p} .
\label{def:norm**}
\end{equation}
The Fenchel\--Young and H\"older inequalities,
\eqref{ineq:Young} and \eqref{ineq:Hoelder}, respectively,
remain valid with ${Q}_V^{**}(u)$ in place of ${Q}_V$.
In particular, we have
\begin{equation}
  | \langle u,f\rangle |\leq \| u\|_{**}\, \| f\|_{*}
  \leq \left( p\, {Q}_V(u) \right)^{1/p}\, \| f\|_{*}
    \quad\mbox{ for } f\in X^{*} ,\
    u\in \core .
\label{ineq:Hoelder**}
\end{equation}
It follows \cite[Ch.~I, \S 4]{Eke-Temam}, that for all $f\in
X^{*}$,
\begin{equation}
    \inf_{ u\in \core }
    \left[{Q}_V(u) - \langle u,f\rangle\right]
  = \inf_{ u\in \core }
    \left[ {Q}^{**}_V(u) - \langle u,f\rangle\right]
  = - {Q}^{*}_V(f) .
\label{eq:inf=inf**}
\end{equation}
\begin{definition}\label{def-gen-min}\nopagebreak
\begingroup\rm
Given a distribution $f\in X^{*}$, we say that a function
$u_0\in W^{1,p}_{\mathrm{loc}}(\Omega)$
is a {\em generalized\/} (or {\em relaxed\/}) {\em minimizer\/}
for the functional
$u\mapsto {Q}_V(u) - \langle u,f\rangle$
if it has the following three properties:
\begin{enumerate}
\renewcommand{\labelenumi}{(\roman{enumi})}
\item[{\rm (i)}]
$\;$
$u_0$ is a (true) minimizer for the functional
\begin{equation*}
  u\mapsto {Q}^{**}_V(u) - \langle u,f\rangle
  \colon X^{**}\to \mathbb{R},
\end{equation*}
hence, $u_0\in X^{**}$.
\item[{\rm (ii)}]
$\;$ $u_0$ satisfies equation $Q'_V(u)=f$ in the sense of
distributions on $\Omega$.
\item[{\rm (iii)}]
$\;$ There exists a (minimizing) sequence $\{ u_k\}_{k=1}^\infty
\subset \core$ such that
\begin{math}
  {Q}_V(u_k) - \langle u_k,f\rangle
  \to - {Q}_V^{*}(f)
\end{math}
and
\begin{equation}
  {Q}_V'(u_k) - f\equiv
  {Q}_V'(u_k) - \langle\,\cdot\, ,f\rangle \to 0
  \quad\mbox{ strongly in } W^{-1,p'}_{\mathrm{loc}}(\Omega)
\label{e:P-S:u_k}
\end{equation}
as $k\to \infty$, together with $u_k\to u_0$ strongly in
$W^{1,p}_{\mathrm{loc}}(\Omega)$, and $u_k
\stackrel{*}{\rightharpoonup} u_0$ weakly\--star in $X^{**}$ as
$k\to \infty$.
\end{enumerate}
\endgroup
\end{definition}
We can now formulate the existence result.
\begin{theorem}[{\cite[Theorem~4.3]{TakTin}}]\label{thm-Minimizer}
Let $\Omega\!\subset \!\R^d$ be a domain, $1\! <\! p\! <
\!\infty$, and $V\in L^{\infty}_{\mathrm{loc}}(\Omega)$. Assume
that the nonnegative functional\/ ${Q}_V$ is weakly coercive.
Then, for every\/ $f\in X^{*}$, the functional\/ $u\mapsto
{Q}_V(u) - \langle u,f\rangle$ is bounded from below on $\core$
and has a generalized minimizer\/ $u_0$ in $X^{**}$ (${}\subset
W^{1,p}_{\mathrm{loc}}(\Omega)$). In particular, this minimizer
verifies the equation $Q'_V(u)=f$.
\end{theorem}

\mysection{Criticality theory}\label{sec:crit} In this section we
discuss positivity properties of the functional $Q$ from
\cite{ky3,ky6} along the lines of criticality theory for
second-order linear elliptic operators \cite{P0,P2}. We note that
Theorem~\ref{thm:VW} and Example~\ref{ex:17} are new results.
\begin{proposition}
\label{monPot} Let  $V_i\in L^\infty_{\mathrm{loc}}(\Omega)$. If
$V_2\gneqq V_1$ and $Q_{V_1}\ge 0$, then $Q_{V_2}$ is subcritical
(weakly coercive).
\end{proposition}

\begin{proposition}
\label{monDom}
Let $\Omega_1\subset\Omega_2$ be domains in $\R^d$ such that
$\Omega_2\setminus\overline{\Omega_1}\neq\emptyset$.
Let $Q_V$ be defined on $C_0^\infty(\Omega_2)$.

1. If $Q_V\ge 0$ on $C_0^\infty(\Omega_2)$, then $Q_V$ is
subcritical in $\Omega_1$.

2. If $Q_V$ is critical in $\Omega_1$, then $Q_V$ is nonpositive
in $\Omega_2$.
\end{proposition}

\begin{proposition}\label{Prop2}
 Let $V_0, V_1\in L^\infty_{\mathrm{loc}}(\Omega)$,
$V_0\neq V_1$. For $t\in \mathbb{R}$ we denote \be
Q_t(u):=tQ_{V_1}(u)+(1-t)Q_{V_0}(u),\ee and suppose that
$Q_{V_i}\geq 0$ on $\core$ for $i=0,1$.

Then $Q_t\geq 0$ on $\core$ for all $t\in[0,1]$. Moreover, $Q_t$
is subcritical in $\Omega$ for all $t\in(0,1)$.
\end{proposition}
\begin{proposition}\label{strictpos}
Let $Q_V$ be a subcritical functional in $\Omega$. Consider
$V_0\in L^\infty(\Omega)$ such that $V_0\ngeq 0$ and $\supp
V_0\Subset\Omega$. Then there exist $\tau_+>0$ and $-\infty\leq
\tau_-<0$ such that $Q_{V+tV_0}$ is subcritical in $\Omega$ for
$t\in(\tau_-,\tau_+)$, and $Q_{V+\tau_+ V_0}$ is critical in
$\Omega$.
\end{proposition}
\begin{proposition}\label{propintcond}
Assume that $Q_V$ admits a ground state $v$ in $\Omega$. Consider
$V_0\in L^\infty(\Omega)$ such that $\supp V_0\Subset\Omega$. Then
there exists $0<\tau_+\leq\infty$ such that $Q_{V+tV_0}$ is
subcritical in $\Omega$ for $t\in(0,\tau_+)$ if and only if
 \be\label{intcond}
  \int_\Omega
V_0|v|^p\dx>0.
 \ee
\end{proposition}
In propositions~\ref{strictpos} and~\ref{propintcond} we assumed
that the perturbation $V_0$ has a compact support. In the
following we consider a wider class of perturbations. Assume that
$Q$ is subcritical in $\Omega$ and $d>p$. It follows from
Theorem~\ref{thmky3} that there exists a continuous weight
function $W$ such that the following Hardy-Sobolev-Maz'ya
inequality is satisfied \be \label{gap_p*1}
Q(u)\ge\left(\int_\Omega W(x)|u(x)|^{p^*}\dx\right)^{p/p^*} \qquad
\forall u\in\core, \ee where $p^*={pd}/{(d-p)}$ is the critical
Sobolev exponent.

The following theorem shows that for a certain class of potentials
$\tilde{V}$ the above Hardy-Sobolev-Maz'ya inequality is preserved
with the same weight function $W$. This theorem extends the
analogous result for $p=2$ proved in \cite[Theorem~2.9]{ky6}. We
may say that such $\tilde{V}$ are {\em small perturbations} of the
functional $Q$ in $\Omega$.
\begin{theorem}
\label{thm:VW} Let $Q$ be the functional \eqref{Q} with $d>p$ and
suppose that \be \label{WW} Q(u)\ge\left(\int_\Omega
W|u|^{p^*}\dx\right)^{p/p^*}  \qquad   u\in C_0^\infty(\Omega),
\ee where $W$ is some positive continuous function (see
Theorem~\ref{thmky3}, and in particular \eqref{gap_p*}). Let \be
\label{VW} \tilde{V}\in L^{\infty}_{\mathrm{loc}}(\Omega)\cap
L^{d/p}(\Omega; W^{-d/p^*}). \ee Consider the one-parameter family
of functionals $\tilde{Q}_\lambda$ defined by
$$
\tilde{Q}_\lambda(u):= Q(u)+\lambda\int_\Omega \tilde{V}|u|^p\dx
\qquad u\in\core,
$$
where $\lambda \in\R$, and let

$$S:= \{\lambda\in\R\mid \tilde{Q}_\lambda\geq 0
 \mbox{ on } \core\}.$$

 (i) If $\lambda \in S$ and $\tilde{Q}_\lambda$ is subcritical in $\Omega$, then there exists
$C>0$ such that  \be \label{subcrt} \tilde{Q}_\lambda(u)\geq
C\left(\int_\Omega W|u|^{p^*}\dx\right)^{p/p^*} \qquad u\in\core.
\ee

(ii) If $\lambda \in S$ and  $\tilde{Q}_\lambda$ admits a ground
state $v$, then for every $\psi\in\core$ such that
$\int_\Omega\psi v\dx\neq 0$ there exist $C, C_1>0$ such that the
following Hardy-Sobolev-Maz'ya-Poincar\'e type inequality holds
\be \label{crtcal} \tilde{Q}_\lambda(u)+C_1\left|\int_\Omega\psi
u\dx\right|^p \ge C\left(\int_\Omega W|u|^{p^*}\dx\right)^{p/p^*}
\qquad u\in\core. \ee

(iii) The set $S$ is a closed interval with a nonempty interior
which is bounded if and only if $\tilde{V}$ changes its sign on a
set of a positive measure in $\Omega$. Moreover, $\lambda \in
\partial S $ if and only if $\tilde{Q}_\lambda$ is critical in $\Omega$.
\end{theorem}
\begin{proof}
(i)--(ii) Assume first that $\tilde{Q}_\lambda$ is subcritical in
$\Omega$, and \eqref{subcrt} does not hold, then there exists a
sequence $\{u_k\}\subset \core$ of nonnegative functions such that
$\tilde{Q}_\lambda(u_k)\to 0$, and $\int_\Omega
W|u_k|^{p^*}\dx=1$. In light of \eqref{gap_max} (with another
weight function $\tilde{W}$) it follows that that $u_k\to 0$ in
$W^{1,p}_{\mathrm{loc}}(\Omega)$.

If $\tilde{Q}_\lambda$ has a ground state $v$, and \eqref{crtcal}
does not hold. Then there exists a sequence $\{u_k\}\subset \core$
of nonnegative functions  such that $\tilde{Q}_\lambda(u_k)\to 0$,
$\int_\Omega\psi u_k\dx\to 0$, but $\int_\Omega
W|u_k|^{p^*}\dx=1$.  It follows from \eqref{Poinc} (with another
weight function $\tilde{W}$) that $u_k\to 0$ in
$W^{1,p}_{\mathrm{loc}}(\Omega)$.

Consequently, for any $K\Subset \Omega$ we have \be\label{intK}
\lim_{k\to\infty}\int_{K} |\tilde{V}||u_k|^p\dx=0. \ee

On the other hand, \eqref{VW} and  H\"older inequality imply that
for any $\varepsilon>0$ there exists $K_\varepsilon\Subset \Omega$
such that
 \be\label{holder}
\left|\int_{\Omega\setminus K_\varepsilon}\!\!\!\!
|\tilde{V}||u_k|^p\dx\right|\!\leq \!\left(\int_{\Omega\setminus
K_\varepsilon}\!\!\!|\tilde{V}|^{d/p}W^{-d/p^*}\dx\right)^{p/d}
\!\!\left(\int_\Omega
\!\!W|u_k|^{p^*}\dx\right)^{p/{p}^*}\!\!\!\!<\varepsilon.
 \ee
Therefore, $\int_\Omega |\tilde{V}||u_k|^p\dx\to 0$.
 Since \be\label{eq17} Q(u_k)\leq
\tilde{Q}_\lambda(u_k)+|\lambda|\int_{\Omega}
|\tilde{V}||u_k|^p\dx,
 \ee
it follows that $Q(u_k)\to 0$. Hence, \eqref{WW} implies that
$\int_\Omega W|u_k|^{p^*}\dx\to 0$, but this contradicts the
assumption $\int_\Omega W|u_k|^{p^*}\dx=1$. Consequently,
\eqref{subcrt}  (resp. \eqref{crtcal}) holds true.

\vskip 3mm

(iii) It follows from Proposition~\ref{Prop2} that $S$ is an
interval, and that $\lambda \in \mathrm{int}\,S$ implies that
$\tilde{Q}_\lambda$ is subcritical in $\Omega$. The claim on the
boundedness of $S$ is trivial and left to the reader.

On the other hand,  suppose that for some $\lambda \in\R$ the
functional  $\tilde{Q}_\lambda$ is subcritical. By part (i),
$\tilde{Q}_\lambda$ satisfies \eqref{subcrt} with weight $W$.
Therefore, \eqref{holder} (with $K_\varepsilon=\emptyset$) implies
that
\begin{equation}\label{0}
 \tilde{Q}_\lambda(u)\geq C\left(\int_\Omega
W|u|^{p^*}\dx\right)^{p/p^*}\geq C_1 \left|\int_\Omega
\tilde{V}|u|^{p}\dx\right| \qquad u\in \core.
\end{equation}
Therefore, $\lambda\in \mathrm{int}\,S$.  Consequently,
$\lambda\in
\partial S$ implies that  $\tilde{Q}_{\lambda}$ is critical in $\Omega$. In particular, $0\in
\mathrm{int}\,S $.

\end{proof}
\begin{example}\label{ex:17}{\em
Let $2\leq p<d$, and let $\Omega\subset\R^d$ be a bounded convex
domain with a smooth boundary. Consider the Hardy functional
$$
Q(u):=\int_\Omega|\nabla
u|^p\dx-\left|\frac{p-1}{p}\right|^p\int_\Omega\dfrac{|u|^p}{d(x,\partial\Omega)^p}\dx.
$$

By \cite{FMT}, the functional $Q$ satisfies the following
Hardy-Sobolev-Maz'ya type inequality
\begin{equation}
Q(u)\ge C\left(\int_\Omega|u|^{p^{*}}\dx \right)^{p/p^*}.
\end{equation}
Let $V\in L^{\infty}_{\mathrm{loc}}(\Omega)\cap L^{d/p}(\Omega)$
be a positive function. By \cite{MS}, there exists a constant
$\lambda_*>0$ such that
$$Q_\lambda(u):= Q(u)-\lambda\int_\Omega V|u|^p\dx$$
is nonnegative for all $\lambda\leq \lambda_*$. Now,
Theorem~\ref{thm:VW} implies that the functional
$Q_{\lambda_*}(u)$ is critical in $\Omega$. Moreover, for
$\lambda<\lambda_*$ the following inequality holds
\begin{eqnarray*}
 Q_\lambda(u)
\geq C_{\lambda}\left(\int_\Omega|u|^{p^*}\dx
\right)^{p/p^*}\qquad u\in \core
\end{eqnarray*}
for some  $C_{\lambda}>0$.

Furthermore,  for every nonzero nonnegative $\psi\in\core$ the
following inequality holds
\begin{eqnarray*}
Q_{\lambda_*}(u)+ C_1\left|\int_\Omega u\psi\dx\right|^p \geq
C\left(\int_\Omega |u|^{p^*}\dx\right)^{p/p^*}\qquad u\in \core
\end{eqnarray*}
for some  $C, C_1>0$.
 }
\end{example}

\mysection{The linear case ($p=2$)}\label{sec:linear case} Some of
the results in the preceding sections have stronger counterparts
in the linear case (for a recent review on the theory of positive
solutions of second-order linear elliptic PDEs, see \cite{P} and
the references therein). In particular, there are several
properties that are true for the linear case but are generally
false or unknown in the general case. This refers in particular to
SCP whose scope of validity when $p\neq 2$ is not completely
understood, and to the convexity of the functional $Q$, which is
known to be generally false for $p>2$ (see references at the
beginning of Section 7).
 On the other hand, as in \cite{Pl,ky2}, we can actually
consider in the linear symmetric case the following somewhat more
general functional than $Q_V$ of the form \eqref{Q}.

Let $A:\Omega \rightarrow \mathbb{R}^{d^2}$ be a measurable
symmetric matrix valued function such that for every compact set
$K\Subset \Omega$ there exists $\mu_K>1$ so that \be \label{stell}
\mu_K^{-1}I_d\le A(x)\le \mu_K I_d \qquad \forall x\in K,
 \ee
where $I_d$ is the $d$-dimensional identity matrix, and the matrix
inequality $A\leq B$ means that $B-A$ is a nonnegative matrix on
$\mathbb{R}^d$. Let $V\in L^{q}_{\mathrm{loc}}(\Omega)$ be a real
potential, where $q>{d}/{2}$. We consider the quadratic form \be
\label{assume}
\mathbf{a}_{A,V}[u]:=\frac{1}{2}\int_\Omega\left(A\nabla u\cdot
\nabla u+V|u|^2\right)\mathrm{d}x   \ee on $\core$
 associated with the Schr\"odinger
equation \be \label{divform}
Pu:=(-\nabla\cdot(A\nabla)+V)u=0\qquad \mbox{ in } \Omega. \ee
 We say that $\mathbf{a}_{A,V}$ is {\em nonnegative} on $C_0^\infty(\Omega)$,
 if $\mathbf{a}_{A,V}[u]\geq 0$ for all $u\in \core$.

Let $v$ be a positive solution of the equation $Pu=0$ in $\Omega$.
Then by \cite[Lemma~2.4]{ky2} we have the following analog of
\eqref{QL1}. For any nonnegative $w\in \core$ we have
 \begin{equation}
\label{e1:2} \mathbf{a}_{A,V}[v w]= \frac{1}{2}\int_\Omega v^2
A\nabla w\cdot\nabla w \,\mathrm{d}x.
 \end{equation}
Moreover, it follows from \cite{ky2,ky3} that all the results
mentioned in this paper concerning the functional $Q$ are also
valid for the form $\mathbf{a}_{A,V}$.
\subsection{Liouville-type theorem}
In the linear case we have the following stronger Liouville-type
statement (cf. Theorem~\ref{thm:main} for $p\neq 2$).
\begin{theorem}[\cite{Pl}]\label{thm:p=2} Let $\Omega$ be a domain in $\R^d$, $d\geq 1$.
Consider two strictly elliptic Schr\"odinger operators with real
coefficients defined on $\Omega$ of the form
\begin{equation}\label{eqpj}
P_j:=-\nabla\cdot(A_j\nabla)+V_j\qquad j=0,1,
\end{equation}
where  $V_j\in L^{p}_{\mathrm{loc}}(\Omega)$ for some $p>{d}/{2}$,
and $A_j:\Omega \rightarrow \mathbb{R}^{d^2}$ are measurable
symmetric matrices satisfying \eqref{stell}.

\vskip 3mm

 \noindent Assume that the following assumptions hold true.
\begin{itemize}
\item[(i)] The operator  $P_1$ admits a ground state
 $\varphi$ in $\Omega$.

\item[(ii)]  $P_0\geq 0$ on $\core$, and there exists a
real function $\psi\in H^1_{\mathrm{loc}}(\Omega)$ such that
$\psi_+\neq 0$, and $P_0\psi \leq 0$ in $\Omega$.

\item[(iii)] The following matrix inequality holds
\begin{equation}\label{psialephia}
(\psi_+)^2(x) A_0(x)\leq C\varphi^2(x) A_1(x)\qquad  \mbox{ a. e.
in } \Omega,
\end{equation}
where $C>0$ is a positive constant.
\end{itemize}
Then the operator $P_0$ is critical in $\Omega$, and $\psi$ is the
corresponding ground state. In particular, $\psi$ is (up to a
multiplicative constant) the unique positive supersolution of the
equation $P_0u=0$ in $\Omega$.
\end{theorem}
%
\subsection{The space $\mathcal D^{1,2}_{A,V}$}
When $p=2$ and the quadratic form $\mathbf{a}_{A,V}$ defined by
\eqref{assume} is nonnegative, $\mathbf{a}_{A,V}$ induces a scalar
product on $\core$. One can regard as the natural domain of the
functional
$$\mathbf{a}_{A,V}^f(u):=\mathbf{a}_{A,V}[u]-\int_\Omega uf\dx$$ a
linear space in which the functional $\mathbf{a}_{A,V}^f$ has a
minimizer for all $f$ such that $\mathbf{a}_{A,V}^f$ is bounded
from below. The functional $\mathbf{a}_{A,V}^f$ is bounded from
below if and only if $\int_\Omega uf\dx$ is a continuous
functional with respect to the norm $(\mathbf{a}_{A,V}[u])^{1/2}$.
The minimum for such $f$ is not attained on $\core$ due to the
strong maximum principle, but any minimizing sequence for
$\mathbf{a}_{A,V}^f$ is a Cauchy sequence. Thus, the natural
domain of $\mathbf{a}_{A,V}$ is the completion of $\core$ in the
norm $(\mathbf{a}_{A,V}[u])^{1/2}$. In the subcritical case, due
to \eqref{gap_max} (which is valid also for subcritical operators
of the form \eqref{divform}), this completion is continuously
imbedded into $W^{1,2}(\Omega; W)$ for some positive continuous
function $W$. By analogy with the classical space $\mathcal
D^{1,2}$, we denote the completion space with respect to the above
norm in the subcritical case by $\mathcal D^{1,2}_{A,V}(\Omega)$
(see \cite{ky2}).

If, however, $\mathbf{a}_{A,V}$ has a ground state $v$, the span
of $v$ becomes obviously the zero element of the completion space
with respect to the norm $(\mathbf{a}_{A,V}[u])^{1/2}$. Recalling
the definition of $\mathcal D^{1,2}(\R^d)$ for $d=1,2$, where the
ground state of $\mathbf{a}_{I,0}[u]=\int_{\R^d}|\nabla u|^2\dx$
is $\mathbf{1}$, and in light of \eqref{Poinc}, we define in the
critical case the norm \be \label{normcr}
\|u\|:=\left(\mathbf{a}_{A,V}[u]+\left|\int_\Omega \psi u\dx
\right|^2\right)^{1/2}, \ee where $\psi$ is any $\core$-function
satisfying $\int_\Omega\psi v\dx\neq 0$. Hence, also in the
critical case the completion of $\core$ with respect to the norm
defined by \eqref{normcr} (which we also denote by  $\mathcal
D^{1,2}_{A,V}(\Omega)$) is continuously imbedded into the function
space $W^{1,2}(\Omega; W)$, with an appropriate weight function
$W$.
\begin{example}
\label{ex:HS}{\em Let $\Omega=\R^d=\R^n\times (\R^m\setminus
\{0\})$, $1\leq m\leq d$, and denote points in $\Omega$ by
$(x,y)\in \R^n\times (\R^m\setminus \{0\})$. Let \be
\label{eq:Hardy}
 \mathbf{a}[u]=\frac12\!\!\int_{\Omega}|\nabla
u|^2\dx\dy\!-\!\left(\frac{m-2}{2}\right)^2\!\!\int_{\Omega}\frac{u^2}{|y|^2}\dx\dy\qquad
u\in C_0^\infty(\Omega). \ee The functional \eqref{eq:Hardy} is
nonnegative due to the Hardy inequality and the constant
$\left(\frac{m-2}{2}\right)^2$ is the maximal constant for which
this is true. Furthermore, if $m<d$, the functional
\eqref{eq:Hardy} is weakly coercive, while for $m=d$ it has a
generalized ground state $v(y)=|y|^{(2-m)/2}$.

 It follows
that completion of $C_0^\infty(\Omega)$ in the norm induced by
\eqref{eq:Hardy} for $m<d$ defines a natural domain for the
functional \eqref{eq:Hardy}. It should be noted, however, that on
the complete space the integrals in the expression
\eqref{eq:Hardy} might be infinite. A more explicit
characterization of the natural domain in this case can be
obtained by noting that $v(x,y)=|y|^{(2-m)/2}$ is a positive
solution of the corresponding equation. Thus, for the functional
\eqref{eq:Hardy} the positive Lagrangian identity \eqref{QL} gives
\be \label{PH} \mathbf{a}[u]=\frac12 \int_{\Omega}|y|^{2-m}|\nabla
(|y|^{(m-2)/2}u)|^2\dx\dy. \ee It follows that the completion
space can be characterized as a space of measurable functions with
measurable weak derivatives for which the integral in \eqref{PH}
is finite. }
\end{example}

An analogous definition of the natural domain for the functional
$Q_V$ could be given for general $p$ whenever the functional $Q_V$
is convex. We should point out that while this is in general
false, one may require the convexity of another functional $\hat
Q$, bounded by $Q_V$ from above and from below. In particular, one
can look at the functionals \eqref{p<2} or \eqref{p>2}. The
following statement from \cite{ky6}) characterizes a convexity
property of $\hat Q$ given by the right hand side of \eqref{p>2}.

\begin{proposition}
\label{convexity} Let $p>2$, and let $v\in
C^1_{\mathrm{loc}}(\Omega)$ be a fixed positive function. Then the
functional
$$\mathcal{Q}(u):=\int_\Omega \left[v^p|\nabla
(u^{2/p})|^p+v^2|\nabla v|^{p-2} u^{2(p-2)/p}|\nabla
(u^{2/p})|^2\right]\dx$$
 is convex on $\{u\in\core, u\ge 0\}$.
\end{proposition}

\subsection{Positive solutions of minimal growth}\label{ssec:mingr2}
We characterize now positive solutions of minimal growth in a
neighborhood of infinity of $\Omega$ in terms of a modified null
sequence of the form $\mathbf{a}_{A,V}$ (cf.
Section~\ref{sec:mingr1}).
\begin{theorem}[{\cite[Theorem~6.1]{ky5}}]
\label{thmmin_null:2}
Suppose that $\mathbf{a}_{A,V}$ is nonnegative on $\core$. Let
$\Omega_1\Subset \Omega$ be an open set, and let $u\in
C(\Omega\setminus \Omega_1)$ be a positive solution of the
equation $Pu=0$ in $\Omega\setminus \overline{\Omega_1}$.
\par
Then $u\in\M_{\Omega,\overline{\Omega_1}}$ if and only if for
every smooth open set $\Omega_2$ satisfying $\Omega_1 \Subset
\Omega_2\Subset \Omega$, and an open set $B\Subset(\Omega\setminus
\overline{\Omega_2})$ there exists a sequence
$\{u_k\}\subset\core$,  $u_k\ge 0$, such that for all $k\in\N$,
$\int_B|u_k|^2\dx=1$, and \be
\lim_{k\to\infty}\int_{\Omega\setminus \overline{\Omega}_2} u^2
A\nabla u_k\cdot\nabla u_k\dx =0.\ee
\end{theorem}
Consider now the following Phragm\'{e}n--Lindel\"{o}f-type
principle that holds in unbounded or nonsmooth domains, and for
irregular potential $V$, provided the subsolution satisfies a
certain decay property (of variational type) in terms of the
Lagrangian $L$ (cf. \cite{Agmon84,LLM,MMP,PS} and
Section~\ref{sec:mingr1}).
\begin{theorem}[Comparison Principle \cite{ky5}]\label{thm_comparison:2}
Assume that $P$ is a nonnegative Schr\-\"odinger operator of the
form \eqref{divform}. Fix smooth open sets $\Omega_1\Subset
\Omega_2\Subset \Omega$. Let $u,v\in
W^{1,p}_{\mathrm{loc}}(\Omega\setminus \Omega_1)\cap
C(\Omega\setminus \Omega_1)$ be, respectively, a positive
subsolution and a supersolution of the equation $Pw=0$ in
$\Omega\setminus \overline{\Omega_1}$ such that $u\leq v$ on
$\partial \Omega_2$.
\par
Assume further that $Pu\in L^\infty_{\mathrm{loc}}(\Omega\setminus
\Omega_1)$, and that there exist an open set
$B\Subset(\Omega\setminus \overline{\Omega_2})$ and a sequence
$\{u_k\}\subset\core$, $u_k\ge 0$, such that  \be
\int_B|u_k|^p\dx=1 \quad \forall k\geq 1, \mbox{ and }\quad
\lim_{k\to\infty}\int_{\Omega\setminus \overline{\Omega}_1} {u}^2
A\nabla ({u}_k/u)\cdot\nabla ({u}_k/u) \,\mathrm{d}x =0.\ee Then
$u\leq v$ on $\Omega\setminus \Omega_2$.
\end{theorem}
\begin{remark}\label{remsub:2}
{\em In Theorem~\ref{thm_comparison:2} we assumed that the
subsolution $u$ is strictly positive. It would be useful to prove
the above comparison principle under the assumption that $u\geq 0$
(cf. \cite{LLM}).}
\end{remark}
\newpage
\begin{center}
{\bf Acknowledgments} \end{center} Y.~P. acknowledges the support
of the Israel Science Foundation (grant No. 587/07) founded by the
Israeli Academy of Sciences and Humanities, by the Fund for the
Promotion of Research at the Technion, and by the Technion
President's Research Fund.

\end{document}